\documentclass[aop,nameyear,MSNbibl,dvips]{arximspdf}
\usepackage{mathbh}

\doi{10.1214/10-AOP610}
\volume{40}
\issue{1}
\pubyear{2012}
\firstpage{130}
\lastpage{145}

\makeatletter

\newcommand{\PP}{\mathbb{P}}
\newcommand{\RR}{\mathbb{R}}
\newcommand{\EE}{\mathbb{E}}
\newcommand{\NN}{\mathbb{N}}

\newcommand{\Var}{\mathsf{Var}}
\newcommand{\II}{\mathbh{1}}

\newtheorem{lemm}[theo]{Lemma}
\newtheorem{prop}[theo]{Proposition}
\newtheorem{theo}{Theorem}[section]
\newtheorem{conj}[theo]{Conjecture}

\newcommand{\eps}{\varepsilon}

\makeatother

\begin{document}
\begin{frontmatter}

\title{Sharp threshold for percolation on expanders}
\runtitle{Sharp threshold for percolation on expanders}

\begin{aug}
\author[A]{\fnms{Itai} \snm{Benjamini}},
\author[B]{\fnms{St\'{e}phane} \snm{Boucheron}},
\author[C]{\fnms{G\'{a}bor} \snm{Lugosi}} and
\author[D]{\fnms{Rapha\"{e}l} \snm{Rossignol}\corref{}\ead[label=e1]{raphael.rossignol@math.u-psud.fr}}

\runauthor{Benjamini, Boucheron, Lugosi and Rossignol}
\affiliation{The Weizmann Institute, Universit\'{e} Paris-Diderot,
Pompeu Fabra University and Universit\'{e} Paris-Sud and CNRS}
\address[A]{I. Benjamini\\
The Weizmann Institute of Science\\
Rehovot POB 76100\\
Israel}
\address[B]{S. Boucheron\\
Laboratoire Probabilit\'{e}s\hspace*{50.32pt}\\
\quad et Mod\`{e}les Al\'{e}atoires\\
Universit\'{e} Paris-Diderot\\
175 rue du Chevaleret\\
75013-F Paris\\
France}
\address[C]{G. Lugosi\\
ICREA and Department of Economics\\
Pompeu Fabra University\\
25-27 Ramon Trias Fargas\\
Barcelona\\
Spain}
\address[D]{R. Rossignol\\
Laboratoire de Math\'{e}matiques\\
UMR 8628\\
Universit\'{e} Paris-Sud\\
Orsay, F91405\\
France\\
\printead{e1}\\
and\\
CNRS, UMR 8628\\
Orsay, F91405\\
France}
\end{aug}

\received{\smonth{7} \syear{2009}}
\revised{\smonth{6} \syear{2010}}

%
\begin{abstract}
We study the appearance of the giant component in random subgraphs
of a given large finite graph $G=(V,E)$ in which each edge is
present independently with probability $p$.
We show that if $G$ is an expander with vertices of bounded degree,
then for any $c\in\ ]0,1[$,
the property that the random subgraph contains a giant component
of size $c|V|$ has a sharp threshold.
\end{abstract}

%

\end{frontmatter}

\section{Introduction}

Percolation theory studies the presence of ``giant''
connected components in a large graph $G$ whose edges are deleted
independently at
random. The meaning of ``giant'' differs whether one looks at
infinite or finite graphs. For the former, it means an infinite
connected component while for finite graphs a giant component usually
means a component of size linear in
the number of vertices of the original graph $G$. Whether the graph is
finite [see \citet{Pittelregular08}, \citet{NachmiasPeres09}] or
infinite [see
\citet{BenjaminiSchramm96beyond}], symmetry
of $G$ has often played a key role in the study of
percolation. In this paper, we are concerned with weakening these
symmetry assumptions in the case of finite graphs, replacing them by
a more geometric assumption. To do so, we follow the path of
\citet{AlonBenjaminiStacey04}, where finite
graphs satisfying an isoperimetric inequality (the so-called
``expanders'') are studied without any symmetry assumption.

As Alon, Benjamini, and Stacey, we study percolation in expanders of bounded
degree.
Consider a finite graph $G_n=(V_n,E_n)$ with $|V_n|=n$ vertices,
and let $G_n(p)$
denote the spanning subgraph of $G_n$ obtained by retaining each edge of
$G_n$ independently with probability $p$. In this paper we consider
``large'' graphs, and the term ``\textit{with high probability}'' refers
to events whose probability is bounded from below by a number
that tends to $1$ as $n\to\infty$.
For any two sets of vertices
$A$ and $B$ in
$G_n$, let $E_n(A,B)$ be the set of all edges with one endpoint
in $A$ and the other in $B$. The
\textit{edge-isoperimetric number} $c(G_n)$ of $G_n$,
also called its \textit{Cheeger constant}, is defined by
\[
\mathop{\min_{A\subset V_n:}}_{0<|A|\leq n/2}\frac
{|E_n(A,A^c)|}{|A|}.
\]
Let $b,d>0$ be constants.
A \textit{$(b,d)$-expander} is a graph
$G_n=(V_n,E_n)$ such that the maximal degree in $G_n$ is not
greater than $d$, and $c(G_n)>b$.\break A~$d$-regular expander is a
$(b,d)$-expander which is
$d$-regular. Theorem 2.1 in \citet{AlonBenjaminiStacey04}
shows that in an expander, there is, with high probability,
never more than one giant component. This allows one to speak about
the giant component. A more precise statement is obtained in
Theorem 2.8 in \citet{AlonBenjaminiStacey04}, showing that, uniformly
over $p$, with probability tending to~$1$, the second largest
component cannot have size larger than $|n|^{\omega(b,d)}$ for some
$\omega(b,d)<1$. Although the same result is conjectured to hold under
less stringent isoperimetric assumptions, it is important to note the
absence of any symmetry assumption.


Arguably the most interesting phenomenon in random graphs is
the emergence of the giant component as $p$ is increased gradually
from $0$ to $1$.
In \citet{AlonBenjaminiStacey04}, Theorem 3.2, the following result is
shown:\looseness=1
\begin{theo}[{[Alon, Benjamini and Stacey (\citeyear{AlonBenjaminiStacey04})]}]
\label{theo:ABS}
Let $d\geq2$ and let $(G_n)_{n\geq0}$ be a sequence of $d$-regular
$(b,d)$-expanders
with girth $(g_n)\rightarrow\infty$.

If $p > 1/(d -1)$, then there exists a $c > 0$ such that
\[
\lim_{n\to\infty} \PP\bigl(G_n(p) \mbox{ contains a component of size
at least
}c|V_n|\bigr) = 1.
\]
If $p < 1/(d-1)$, then for any $c > 0$,
\[
\lim_{n\to\infty}
\PP\bigl(G_n(p) \mbox{ contains a component of size at least
}c|V_n|\bigr) = 0.
\]
\end{theo}


It is tempting to conjecture that the regularity and high-girth
assumptions in Theorem \ref{theo:ABS} can be removed at the price of
losing the precise location of the threshold, thus showing that
\textit{a giant component emerges in an interval of length $o(1)$ in any
expander}.
\begin{conj}
\label{conj:birth}
Let $G_n$ be a $(b,d)$-expander. There exist $0<q_1<q_2<1$,
depending on $b$ and $d$ and $p^*_n\in[q_1,q_2]$,
such that, for every $\eps>0$,
if $p_n \geq p^*_n+\eps$, then there exists a $c > 0$ such that,
with high probability,
\[
G_n(p_n) \mbox{ contains a component of size at least } cn,
\]
and if $p_n \leq p^*_n-\eps$, then for any $c > 0$,
with high probability,
\[
G_n(p_n) \mbox{ does not contain a component of size at least
}cn.\vadjust{\goodbreak}
\]
\end{conj}

However, it is not entirely clear whether Conjecture \ref{conj:birth}
can hold
without a minimum of homogeneity of the underlying graph $G_n$.
\citet{BenjaminiNachmiasPeres09}, Theorem 1.3,
establishes Conjecture \ref{conj:birth}
under the additional ``homogeneity'' assumption
of weak convergence of the sequence of graphs
$G_n$ to an infinite bounded-degree graph.

The main result of this paper is a sharp threshold result for the
events ``$G_n(p)$ contains a component of order at least $cn$''
for every $c\in\ ]0,1[$. In a sense, it can be seen as a weakening of
Conjecture \ref{conj:birth}. The weakness is that we cannot assert
the existence of the ``threshold function'' $p^*_n$, but we can
still interpret this result as the fact that
\textit{in any expander, every giant component of given proportion
emerges in
an interval of length $o(1)$}. This is formalized as follows.
\begin{theo}
\label{theo:main}
Let $G_n$ be a $(b,d)$-expander and let $c\in\ ]0,1[$.
There exist constants $q_1=q_1(d) >0$ and $q_2=q_2(c) \in\ ]q_1,1[$
and $p^*_n(c)\in[q_1,q_2]$ such that, for every $\eps>0$.

If $p_n \geq p^*_n(c)+\eps$, then, with high probability,
\[
G_n(p_n) \mbox{ contains a component of size at least
}cn.
\]
If $p_n \leq p^*_n(c)-\eps$, then, with high probability,
\[
G_n(p_n) \mbox{ does not contain a component of size at least
}cn.
\]
\end{theo}

The rough idea of the proof is the following. First we show that the
expansion and bounded-degree properties are sufficient to imply that
the value of $p$ for which the probability that $G_n(p)$ contains
a component of size $cn$ equals some fixed constant $\alpha\in(0,1)$
stays bounded away from zero and one.
This may be proved by putting together some arguments
of \citet{BenjaminiSchramm96beyond} and
\citet{AlonBenjaminiStacey04}.
Now to show that the threshold is sharp, it suffices to prove
that the threshold width is bounded by a function of $n$
that tends to zero.
The proof of this fact has two main components. First we show
that around the critical probability, the derivative (with respect to $p$)
of the expected size of the largest component is proportional to
$n$. This is done by using Russo's lemma and the expansion
property. In other words, the expected size of the largest
component grows quickly (with $p$) around the threshold value.
Then to show that the probability of existence of a component of
size $cn$ has a rapid growth around the threshold,
it suffices to show that the size of the largest component
is concentrated in the sense that its standard deviation is of
a smaller order of magnitude than its mean. This may be achieved by
applying a general bound for the variance of a function of independent
Bernoulli variables due to \citet{FalikSamorodnitsky07}.

Note that if the underlying graph $G_n$ is transitive, then
after showing that the critical probability is bounded away
from
zero and one, the theorem follows immediately from a general\vadjust{\goodbreak}
result of \citet{FriedgutKalai} which implies that the
threshold width is at most $O(1/\log n)$. However, because of
the lack of symmetry assumption, the theorem of Friedgut and Kalai
is not directly applicable. Some of the main ingredients of our
proof are similar to those of Friedgut and Kalai as in both proofs
Russo's lemma and an appropriate variance bound (based on either
hypercontractivity or a~logarithmic Sobolev inequality) play a key role.
We use an elegant inequality of Falik and Samorodnitsky that can
also be used to prove the Friedgut--Kalai theorem.

The paper is organized as follows. In Section \ref{sec:notations} we
introduce some notation. In Section \ref{sec:away} we collect some
arguments existing in the literature to show that the thresholds
$p^*_n(c)$ are bounded away from zero and
one. Finally, Section~\ref{sec:main} is devoted to the proof of our main
result, Theorem \ref{theo:main}. In fact, we prove something more in
Theorem \ref{theo:sharpthreshold}, namely that the
threshold width is at most of order $O((\log n)^{-1/3})$.

We end this Introduction with an open question. It would be interesting
to relax the isoperimetric condition and replace the positivity of
the Cheeger constant by a weaker condition of the type
\[
\mathop{\min_{A\subset V_n:}}_{0<|A|\leq n/2}\frac
{|E(A,A^c)|}{|A|^\alpha}>0
\]
for some $\alpha<1$.

\section{Notation}
\label{sec:notations}

Let $b>0$ and $d$ a positive integer.
Throughout the paper, $G_n=(V_n,E_n)$
denotes a $(b,d)$-expander.
For a subset $W$ of vertices, we denote by
$\partial_EW = E(W,W^c)$ the exterior edge-boundary of $W$.

Each point (or ``configuration'')
$x \in\{0,1\}^{E_n}$ is identified with the subgraph of
$G_n$ with vertex set $V_n$ and edge set obtained by removing
from~$E_n$ all edges $e$ such that $x(e)=0$. For a
given $p\in[0,1]$, we
equip the space $\{0,1\}^{E_n}$ with the product probability measure
$\mu_{n,p}$ under which every $x(e)$ is independently $1$ (resp., $0$)
with probability $p$ [resp., $(1-p)$]. For any function
$f\dvtx\{0,1\}^{E_n}\to\RR$, we denote by $\EE_{n,p}(f)=\int f(x) \,d\mu
_{n,p}(x)$
the mean, and for any $\alpha\geq1$, the norm
\mbox{$\|\cdot\|_{\alpha,p}$} denotes
\[
\|f\|_{\alpha,p}=(\EE_{n,p}(|f|^\alpha))^{
{1}/{\alpha}} .
\]
For $x\in\{0,1\}^{E_n}$ and $i\in\NN^*$, define
$\mathcal{C}^{(i)}_n=\mathcal{C}^{(i)}_n(x)$ to be the
$i$th largest connected component\vspace*{1pt} in the configuration $x$, and let
$L_n^{(i)}=L_n^{(i)}(x)$ be the number of vertices in
$\mathcal{C}^{(i)}_n$. We also denote by $C(v)$ the connected
component containing a vertex $v\in V_n$.

For any $c\in\ ]0,1[$, the subset of $\{0,1\}^{E_n}$ defined by
$L_n^{(1)}\geq c n$ is monotone (i.e., if the inequality
holds for a graph then by adding any new edge it still holds),
and therefore $\mu_{n,p}\{L_n^{(1)}\geq c n\}$ is a strictly increasing
polynomial of~$p$. Thus, for any $\alpha\in[0,1]$, we may
define\vadjust{\goodbreak}
$p_{n,\alpha}(c)$ as the
unique real number~$p$ in $[0,1]$ such that
\[
\mu_{n,p}\bigl\{L_n^{(1)}\geq c n\bigr\}=\alpha.
\]
The threshold function is defined as
\[
p_n^*(c) =p_{n,1/2}(c) .
\]
%
When $n$ is clear from the context, we omit the subscript $n$.

\section{Thresholds of giant components are bounded away from zero and one}
\label{sec:away}

The fact that the values $p_{n,\alpha}(c)$ are bounded
away from zero and one may be proved by putting together
arguments of \citet{BenjaminiSchramm96beyond} and
\citet{AlonBenjaminiStacey04}.
Here we give a self-contained proof. We also add an estimate
for the decay of the probability of not having a~giant component of
fixed size for $p$ close enough to $1$, which seems to us interesting
on its own.
\begin{prop}
\label{prop:criticalpaway01}
Let $c\in\ ]0,1[$. There is
a constant $q_1$, depending only on $d$ and there
exists $q_2(c)\in\ ]q_1,1[$ such that for any
$\alpha\in\ ]0,1[$, for all $n$ large enough,
$p_{n,\alpha}(c)\in\ ]q_1,q_2(c)[$.

Furthermore, for any $c\in\ ]0,1[$, there are strictly positive
constants $C_1$ and~$C_2$,
depending only on $b$ and $d$, such
that for every $p\geq q_2(c)$,
\[
\mu_{n,p}\bigl(L_n^{(1)}\geq cn\bigr)\geq1-C_1e^{-C_2n} .
\]
\end{prop}
\begin{pf}
The fact that $p_{n,\alpha}(c)$ is bounded away from $0$ may be proved using
standard branching-process arguments as follows. Fix
$q_1<1/(d-1)$ and suppose that $p\leq q_1$. Since the degrees are
bounded by $d$, the connected component $C(v)$ of a vertex $v\in V_n$
has a size not larger than $S$, where $S$ is the total number of
descendants of the root in a (sub-critical) Galton--Watson process with
offspring distribution $\mathcal{B}(d-1,p)$ [except for the first
offspring, which has distribution $\mathcal{B}(d,p)$]. Since the binomial
distribution $\mathcal{B}(d-1,p)$ possesses exponential moments, it is
well known [see, e.g, Exercise~5.22 in \citet{LyonsPeresinprogress}] that
there are some values $\lambda>0$, $M<\infty$, depending only on $d$
and $q_1$, such that, for every $n$ and $p\leq q_1$,
\[
\EE_p(e^{\lambda S})\leq M .
\]
Thus, for any $t>0$ and $p\leq q_1$,
%
%
\begin{equation}
\label{eq:subcritical}
\mu_p\bigl(L_n^{(1)}>t\bigr)\leq nMe^{-\lambda t} .
\end{equation}
In particular,
\[
\mu_p\biggl(L_n^{(1)}>\frac{2}{\lambda}\log(nM^{1/2})\biggr)\leq
\frac{1}{n} .
\]
Thus, for any $\alpha\in\ ]0,1[$, for $n$ large enough, $p_{n,\alpha
}(c)> q_1$.

The fact that $p_{n,\alpha}(c)$ is bounded away from $1$ follows
essentially from Theorem~2 and Remark 2
in \citet{BenjaminiSchramm96beyond}\vadjust{\goodbreak} and Lemma~2.2 and
the proof of Proposition 3.1 in \citet{AlonBenjaminiStacey04}.
To detail the proof,
first we show that there are constants $p_0(b)<1$, $a(b)>0$
and $C(b,d)>0$ such that for every $p\geq p_0$ and every $n\in\NN^*$,
%
%
\begin{equation}
\label{eq:ABSprop31bis}
\mu_p\bigl(L_n^{(1)}\geq an\bigr)\geq1-e^{-Cn} .
\end{equation}
Then, we slightly extend Lemma 2.2 in \citet{AlonBenjaminiStacey04} in
proving that for every $c_1\in\ ]0,1/2[$ and $c_2\in\ ]1/2,1[$, there is a
constant $q_3(c_1,c_2)$, depending only on $c_1$, $c_2$, $b$ and
$d$, such that, for every $p>q_3(c_1,c_2)$,
%
%
\begin{eqnarray}
\label{eq:ABSlemme22bis}
&&\PP\bigl(G_n(p)\mbox{ contains a component of size in
}[c_1n,c_2n[\bigr) \nonumber\\[-8pt]\\[-8pt]
&&\qquad\leq4\biggl(1+\frac{1}{c_1}\biggr)e^{-n} .\nonumber
\end{eqnarray}
The proof of this extension follows from an argument similar to
\citet{AlonBenjaminiStacey04}.

Once inequalities (\ref{eq:ABSprop31bis}) and (\ref{eq:ABSlemme22bis})
are proved, Proposition \ref{prop:criticalpaway01} follows
with the choice $q_2(c)=\max\{q_3(\min\{1/4,a\},\max\{3/4,c\}
),p_0(b)\}$.\vspace*{8pt}

\textsc{Proof of (\ref{eq:ABSprop31bis})}.\quad
First we show that if $p>\frac{1}{1+b}$,
then there is some $\delta>0$, depending only on $p-\frac{1}{1+b}$,
such that
%
%
\begin{equation}
\label{eq:BSbis}
\forall v\in V_n\qquad \mu_p\bigl(|C(v)|\geq n/2\bigr)\geq\delta.
\end{equation}
To see this, we construct recursively the component $C(v)$ and its
edge-boundary $W(v)$ as follows. First, we order the edges in
$E_n$. Let $C_1=\{v\}$ and $W_1=\varnothing$. At each step of the
algorithm, $C_k$ denotes a subset of $C(v)$ and~$W_k$ a subset of
$W(v)$. At step $k$, we explore the first (in the aforementioned
order) edge $e_k=(y,z)$ from $E_n\setminus W_k$ that is adjacent to a
vertex $y$ from $C_k$ and to a vertex $z$ from $V_n \setminus C_k$, if
there exists such an edge. Otherwise ($\partial_EC_{k-1}\subset
W_{k-1}$), $C_{k+1}=C_k$ and $W_{k+1}=W_k$. If the edge $e_k=(y,z)$ is
open [i.e., $x(e_k)=1$], let $C_{k+1}=C_{k}\cup\{z\}$ and
$W_{k+1}=W_{k}$. If the edge $(y,z)$ is closed [i.e., $x(e_k)=0$],
let $C_{k+1}=C_{k}$ and $W_{k+1}=W_{k}\cup\{e_k\}$. We have
$C(v)=\bigcup_{k=1}^{\infty}C_k$. If $|C(v)|<n$, then there is a smallest
$N$ such that $W_N=\partial_EC_N$. Notice that $N=|C_N|+|W_N|$ and
$C_N=C(v)$. Thus, by the expansion assumption, if
$|C_N|\leq\frac{n}{2}$,
\[
N - |W_N| = |C_N| \le\frac{|W_N|}{b},
\]
so $|W_N|\ge Nb/(1+b)$.

It is easy to see that, under $\mu_{n,p}$, $(x(e_1),\ldots,x(e_N))$
can be completed so as to form an infinite i.i.d. sequence of Bernoulli
random variables with parameter $p$. Thus, to construct $C(v)$,
we flipped $N-1$ independent $(p,\break 1-p)$-coins
and at least $Nb/(1+b)$ among them turned out zero. But if\break
$p>\frac{1}{1+b}$, then, with positive probability $\delta$,\vadjust{\goodbreak}
depending only on $p-\frac{1}{1+b}$, a~random infinite sequence of
i.i.d. Bernoulli variables of parameter $p$ does not have an $n$ such
that at least $(n+1)b/(1+b)$ among the $N$ first coordinates equal
$0$. The last fact is a consequence of the law of large numbers. This
proves inequality (\ref{eq:BSbis}).


Now, fix $q\in\ ](1+b)^{-1},1[$, let $R$ be a positive real number to
be chosen later and define $S_n$ to be the number of
vertices which belong to a component of size at least
$R/2$
\[
S_n=\sum_{v\in V_n}\II_{|C(v)|>R/2} .
\]
Denoting $X_v=\II_{|C(v)|>R/2 }$, notice that $X_v$ and $X_{v'}$ are
independent as soon
as $d(v,v')>R$
where $d(v,v')$ is the distance of vertices $v$ and $v'$
according to the shortest path metric in $G_n$.
Thus, using the fact that the maximal
degree in $G_n$ is at most $d$,
the maximal degree in the dependency graph of
$(X_v)_{v\in V_n}$ is less than $d^R$.
Recall that the dependency graph of the random variables
$(X_v)_{v\in V_n}$ is given by the vertex set $V_n$ and
the edge set satisfying that if for two disjoint sets of vertices
$A$ and $B$ there is no edge between $A$ and $B$, then
the families $(X_v)_{v\in A}$ and $(X_v)_{v\in B}$ are independent.
Thus, by Theorem 2.1 in \citet{Janson04} for any $t>0$,
\[
\mu_p\bigl(S_n<\EE_p(S_n)-t\bigr)\leq e^{-{2t^2}/({nd^{R}})} .
\]
Notice that a similar result would be obtained from the method of
bounded differences. From (\ref{eq:BSbis}), we see that if $p\geq q$,
$\EE_p(S_n)\geq\delta n$. Choosing $t=\EE_p(S_n)/2$ in the
above inequality gives, for any $p\geq q$,
\[
\mu_p(S_n<\delta n/2)\leq e^{-{\delta^2n}/({2d^{R}})} .
\]
This means that with probability at least $1-e^{-{\delta^2n}/({2d^{R}})}$,
there are at least $\delta n/2$ vertices which belong
to components of size at least $\frac{R}{2} $. Then,
fix $p_0\in\ ]q,1[$. The proof of Proposition 3.1 in
\citet{AlonBenjaminiStacey04} shows that if $R$ is chosen large enough,
there is some positive constant $C$ depending on $b$, $d$, $\delta$,
$p_0$ and
$q$ such that with probability at least $1-e^{-{\delta^2n}/({2d^{R}})}-e^{-C n}$, for $n$ large enough, there is a
component of size at least $\delta n/6$ in $G_n(p_0)$.
We recall their argument:
fix a set of at most $r=\delta n/R$
components of size at least $R/2$ which
contain together at least $\delta n/2$ vertices.
If $\eps=1-\frac{1-p_0}{1-q}$, $G(p_0)$ has the same law as
$G(q)\cup G(\eps)$, where $G(q)$ and $G(\eps)$ are independent.
Then, we claim that there is some $C$ depending on $b$, $d$ and
$\eps$ such that with probability
at least $1-e^{-Can}$, in the random graph~$G(\eps)$, there is no way
of splitting these components into two parts $A$ and $B$, each containing
at least $\delta n/6$ vertices,
with no path of $G(\eps)$ connecting the two parts.
This will imply that, with
the required probability, $G(q)\cup G(\eps)$ contains
a connected component consisting of at least $\delta n/6$
vertices. Now, we show the claim. Let us fix two parts $A$ and $B$ of
the components above, each containing
at least $\delta n/6$ vertices. Thanks to Menger's theorem, there are
at least $b\delta n/6$ edge-disjoint paths between $A$ and $B$. Since
the total number of edges is less than $dn/2$, at least half of these
paths are of length not larger than $6d/(b\delta)$. Thus, the
probability that there is no path between~$A$ and $B$ in $G(\eps)$ is
at most
\[
\bigl(1-\eps^{6d/(b\delta)}\bigr)^{b\delta n/12}\leq e^{-b\delta
n\eps^{6d/(b\delta)}/12} .
\]
Now, there are at most $2^r=2^{\delta n/R}$ ways to choose $A$ and
$B$. Thus, the probability that there is a way to split the components
into two parts~$A$ and~$B$, each containing
at least $\delta n/6$ vertices,
with no path of $G(\eps)$ connecting the two parts is at most
\[
2^{\delta n/R}e^{-b\delta n\eps^{6d/(b\delta)}/12}\leq e^{-b\delta
n\eps^{6d/(b\delta)}/24} ,
\]
as soon as $R$ is larger than $24\log2/(b\eps^{6d/(b\delta)})$. This
finishes the proof of the claim, and thus the one of (\ref{eq:ABSprop31bis}).
\vspace*{8pt}

\textsc{Proof of (\ref{eq:ABSlemme22bis})}.\quad It is well known [see
\citet{FlajoletSedgewick09}, Example~I.14, page 68] that an
infinite $d$-regular rooted tree contains precisely
$\frac{1}{(d-1)r+1}{r\choose dr}$ rooted subtrees of size r. This
number is at most $(de)^r$. To a~graph $G$ of maximum degree less than
$d$ and a
vertex $v$ of that graph, one may associate a subtree of the infinite
$d$-regular tree rooted at $v$ by considering the self-avoiding paths
issued from $v$ in $G$. Through this correspondence, any
connected component of size $r$ of $G$ containing $v$ is mapped to a
different subtree of size $r$ (through the choice of a spanning
tree of the component in~$G$). Thus, the total
number of connected subsets of size~$r$ in~$V_n$ is at most
$n(de)^r/r$. Now, using the expanding property, for any subset $U$
of size $r$, the probability that all edges in $E(U,U^c)$ are
absent is at most $(1-p)^{br}$ if $r\leq n/2$ and at most
$(1-p)^{b(n-r)}$ if $r> n/2$. Thus, the probability that there
is a connected component of size in $[c_1n,c_2n[$ is at most
\begin{eqnarray*}
&&\sum_{r=\lceil c_1n\rceil}^{\lfloor n/2\rfloor}\frac{n(de)^r}{r}
(1-p)^{br}+\sum_{r=\lfloor
n/2\rfloor+1}^{\lfloor c_2n\rfloor}\frac{n(de)^r}{r}(1-p)^{b(n-r)}\\
&&\qquad\leq\frac{1}{c_1}\frac{(de(1-p)^b
)^{c_1n}}{1-(de(1-p)^b)}+2(1-p)^{nb}\frac{(de(1-p)^{-b}
)^{c_2n+1}}{de(1-p)^{-b}-1} \\
&&\qquad\leq
\frac{4}{c}e^{-n}+4e^{-n} ,
\end{eqnarray*}
provided that
\[
de(1-p)^{-b}\geq2,\qquad (de)^{c_2}(1-p)^{b(1-c_2)}\leq\frac{1}{e}
\quad\mbox{and}\quad \bigl(de(1-p)^{b}\bigr)^{c_1}\leq\frac{1}{e} .
\]
These conditions are satisfied if $p$ is larger than some
$q_3(c_1,c_2)<1$. This defines the value of $q_3(c_1,c_2)$ for which
(\ref{eq:ABSlemme22bis}) is valid for every $p>q_3(c_1,c_2)$.
\end{pf}

\section{Threshold phenomenon for the appearance of a giant component}
\label{sec:main}

In this section, we prove our main result, Theorem
\ref{theo:main}. The main step is stated below in Theorem
\ref{theo:sharpthreshold}, showing that the threshold for having a
component of size at least $cn$ has width of order at most
$O((\log n)^{-1/3})$. Theorem \ref{theo:sharpthreshold},
together with Proposition \ref{prop:criticalpaway01} imply our main result,
Theorem \ref{theo:main}. [Recall that $p_n^*(c)=p_{n,1/2}(c)$.]
\begin{theo}
\label{theo:sharpthreshold}
Let $\alpha<1/2$ and $c\in\ ]0,1[$.
There is a constant $C_3$, depending only on $c$,
$\alpha$, $b$, and $d$, such that, for any $n$,
\[
p_{n,(1-\alpha)}(c)-p_{n,\alpha}(c)\leq\frac{C_3}{(\log n)^{1/3}} .
\]
\end{theo}

Here is the idea of the proof. Let us call informally, the
``super-critical phase'' the set of values of $p$ such that a giant
component of size $cn$ has appeared with probability greater than
some $\alpha>0$. The main idea of the proof is to show that for most
values of $p$ in the super-critical phase, the standard deviation of
$L_n^{(1)}$, the size of the largest component, is small with respect
to its mean (which is of the order $n$). This is shown essentially
in Lemma~\ref{lemm:variancederivee} below. In this lemma, we crucially
use an estimate of \citet{AlonBenjaminiStacey04} for the probability
that the second largest component is ``large'' (greater than some
$n^\omega$ with $\omega\in\ ]0,1[$). Next, it follows immediately
from the expanding
property that the mean of $L_n^{(1)}$ has
a~derivative at least of order $n$ inside the super-critical
phase. This is proved in Lemma \ref{lemm:grandederivee}. These two
facts imply that the threshold is sharp: when $p$ goes from
$p_{n,1-\alpha}$ to $p_{n,\alpha}$, the size of the largest component
increases by
a~positive fraction of the number of vertices, since the
fluctuations of this size around its mean is small with respect to
the number of vertices.

Now we turn to the proof which relies on a series of lemmas.
First we need some technical definitions. For any function
$f\dvtx\{0,1\}^{E_n}\to\RR$ and any $e\in E_n$, define the operator
$\Delta_{e,p}$ as
\[
\Delta_{e,p}f(x)=f(x)-\int f(x) \,dx(e) ,
\]
where the integration with respect to $x(e)$ is understood with
respect to the Bernoulli measure with parameter $p$. When there is no
ambiguity, we write~$\Delta_e$ instead of $\Delta_{e,p}$. Finally, the
following notation will be useful: when~$X$ and $X'\in\{0,1\}^{E_n}$
are independent and distributed according to~$\mu_{n,p}$, and $e$
belongs to $E_n$, we denote by $X^{(e)}$ the random
configuration obtained from~$X$ by replacing $X_e$ by $X'_e$.
\begin{lemm}
\label{lemm:norm1}
There exist $\delta(b,d) <1$ and $K(b,d)<\infty$ such that,
for every~$n$,
\[
\sup_{p}\sup_{e\in E_n}\bigl\|\Delta_{e,p}L_n^{(1)}\bigr\|_{1,p}\leq
K(b,d)n^{\delta(b,d)} .
\]
\end{lemm}
\begin{pf}
To lighten notation, we write $f(x)=L_n^{(1)}(x)$ for the size of~the largest component.
A look at the proof of Theorem 2.8 in \citet{AlonBenjaminiStacey04}
reveals that there are three positive real numbers $K_1(b,d)<\infty$,
$\omega(b,d)<1$, and $g(b,d)>0$, depending only on $b$ and $d$, such
that\looseness=-1
\[
\sup_p\PP\pmatrix{G_n(p)\mbox{ contains more than one
component}\cr
\mbox{of size at least }n^{\omega(b,d)}}
\leq K_1(b,d)n^{-g(b,d)} .
\]\looseness=0
The estimate for $p$ close to $0$ is not made explicit in
\citet{AlonBenjaminiStacey04} but follows
from (\ref{eq:subcritical}). Since $f$ is monotone, one may write for
any fixed $p$,
\[
\|\Delta_{e,p}f\|_{1,p} = 2\EE\bigl[\bigl(f(X)-f\bigl(X^{(e)}\bigr)\bigr)_-\bigr]
\]
[where $y_-=\max(0,-y)$ denotes the negative part of a real number $y$].
Let~$A_n$ be the event that the size of the second largest connected
component is at most $n^{\omega(b,d)}$. Notice that
$(f(X)-f(X^{(e)}))_-$ can only be positive if $x(e)=0$, and
the edge $e$ is adjacent to the largest component in $G_n$.
Then the difference between $f(X^{(e)})$ and $f(X)$ is the size of
the component that gets attached to the largest component by
adding the edge $e$ to $G_n$.
Thus $(f(X)-f(X^{(e)}))_-$ is always bounded by
the size of the second largest connected
component in $G_n$. It is also smaller than $n$. Thus, we have
\begin{eqnarray*}
\bigl\|\Delta_{e,p}L_n^{(1)}\bigr\|_{1,p}
&\leq&2\EE\bigl[\bigl(L_n^{(1)}(X)-L_n^{(1)}\bigl(X^{(e)}\bigr)\bigr)_-\bigr] \\[-2pt]
&\leq&2n^{\omega(b,d)}+2n\PP(A_n^c)\\[-2pt]
&\leq&\bigl(2+K_1(b,d)\bigr)n^{\delta(b,d)} ,
\end{eqnarray*}
where $\delta(b,d)=\max\{\omega(b,d),1-g(b,d)\}<1$.\vspace*{-2pt}
\end{pf}

Next we establish an upper bound for the variance of the second largest
component. Our main tool is a result of
\citet{FalikSamorodnitsky07} that gives an improved estimate over
the Efron--Stein inequality for functions defined on the binary hypercube.
Recall that the Efron--Stein inequality implies that if
$f\dvtx\{0,1\}^{E_n} \to\RR$, then
\[
\Var(f) \le\sum_{e\in E_n}\|\Delta_e f\|_{2,p}^2
\]
[see \citet{EfSt81}]. The next inequality appears in this form in
\citet{BenjaminiRossignol07}:\vspace*{-2pt}
\begin{lemm}[(Falik and Samorodnitsky)]
\label{lemmFS}
Let $f$ belong to $L^1(\{0,1\}^{E_n})$. Suppose that
$\mathcal{E}_1(f)$ and $\mathcal{E}_2(f)$ are two real numbers such that
\begin{eqnarray*}
\mathcal{E}_2(f)&\geq&\sum_{e\in E_n}\|\Delta_e f\|_2^2 ,
\\[-2pt]
\mathcal{E}_1(f)&\geq&\sum_{e\in E_n}\|\Delta_e f\|_1^2\vadjust{\goodbreak}
\end{eqnarray*}
and
\[
\frac{\mathcal{E}_2(f)}{\mathcal{E}_1(f)}\geq e .
\]
Then
\[
\Var(f)\leq
2\frac{\mathcal{E}_2(f)}{\log({\mathcal{E}_2(f)}/({\mathcal
{E}_1(f)\log({\mathcal{E}_2(f)}/{\mathcal{E}_1(f)})}))} .
\]
\end{lemm}

This\vspace*{1pt} inequality may be used to derive our next key lemma,
implying that for most values of $p$ in
the super-critical phase, the standard deviation of~$L_n^{(1)}$ is small with respect to its mean.
\begin{lemm}
\label{lemm:variancederivee}
There is a constant $C(b,d)<\infty$ such that, for any $p$ and~$n$,
\[
\Var_p\bigl(L_n^{(1)}\bigr)\leq C(b,d)\frac{n}{\log n}\,\frac{d\EE
_p(L_n^{(1)})}{dp} .
\]
\end{lemm}
\begin{pf}
Fix $\beta\in\ ]0,1-\delta(b,d)[$, where $\delta(b,d)$ is
given in Lemma \ref{lemm:norm1}. Define
\[
\mathcal{E}_2\bigl(L_n^{(1)}\bigr)=\sum_{e\in E_n}\bigl\|\Delta_e
L_n^{(1)}\bigr\|_2^2
\]
and
\[
\mathcal{E}_1\bigl(L_n^{(1)}\bigr)=\sum_{e\in E_n}\bigl\|\Delta_e
L_n^{(1)}\bigr\|_1^2 .
\]
We distinguish two cases depending on the relationship
between $\mathcal{E}_1(L_n^{(1)})$ and $\mathcal{E}_2(L_n^{(1)})$:
\begin{itemize}
\item If $\mathcal{E}_2(L_n^{(1)})\leq n^\beta\mathcal{E}_1(L_n^{(1)})$,
then using
Lemma \ref{lemm:norm1},
\[
\mathcal{E}_2\bigl(L_n^{(1)}\bigr)\leq K(b,d)n^{\beta+\delta(b,d)}
\sum_{e\in E_n}\bigl\|\Delta_e L_n^{(1)}\bigr\|_1 .
\]
Since $L_n^{(1)}$ is a monotone increasing function on the binary hypercube,
a straightforward generalization of
Russo's lemma [see \citet{Rossignol06}] implies that
%
%
\begin{equation}
\label{eq:Russo}
\frac{d\EE_p(L_n^{(1)})}{dp}=\frac{1}{2p(1-p)}\sum_{e\in
E_n}\bigl\|\Delta_e L_n^{(1)}\bigr\|_1 .
\end{equation}
On the other hand, by the Efron--Stein inequality,
\[
\Var_p\bigl(L_n^{(1)}\bigr)\leq\mathcal{E}_2\bigl(L_n^{(1)}\bigr) .
\]
Thus, ignoring the term $2p(1-p) < 1$,
\[
\Var_p\bigl(L_n^{(1)}\bigr)\leq K(b,d)
n^{\beta+\delta(b,d)}\times\frac{d\EE_p(L_n^{(1)})}{dp} .
\]
\item
If $\mathcal{E}_2(L_n^{(1)})> n^\beta\mathcal{E}_1(L_n^{(1)})$, then Lemma
\ref{lemmFS} implies that there is a constant~$C_1$, depending only on
$\beta$ (i.e., on $b$ and $d$), such that
\[
\Var_p\bigl(L_n^{(1)}\bigr)\leq C_1\frac{\mathcal{E}_2(L_n^{(1)})}{\log n} .
\]
But since $L_n^{(1)}$ is positive and always smaller than $n$,
we have
$\mathcal{E}_2(L_n^{(1)})\leq n\sum_{e\in E_n}\|\Delta_e
L_n^{(1)}\|_1$.
Thus, Russo's lemma implies that
\[
\Var_p\bigl(L_n^{(1)}\bigr)\leq C_1\frac{n}{\log n}\,\frac{d\EE
_p(L_n^{(1)})}{dp} .
\]
\end{itemize}
In both cases, the result follows.
\end{pf}

The following easy lemma states that, whatever $\gamma<1$ and
$\eps\in\ ]0,1[$ are, there is always some $c<1$ such that the
probability of having a component of size $cn$ is less than
$\eps$ if $p\leq\gamma$.
\begin{lemm}
\label{lemm:proche1}
Let $\gamma<1$ and $\eps\in\ ]0,1[$. Then there is some $c<1$ such
that, for $n$ large enough,
\[
\mu_{n,\gamma}\bigl(L_n^{(1)}\geq cn\bigr)\leq\eps.
\]
\end{lemm}
\begin{pf}
The size of the largest connected component is less
than $n-N$, where $N$ is the number of isolated vertices (except when all
vertices are isolated, in which case $L_n^{(1)}$ is 1). Let $X_v$
denote the indicator function of the event ``$v$ is isolated.'' If
$d(v,v')\geq2$, $X_v$ and $X_v'$ are independent. Thus, the maximal
degree in the dependency graph of
$(X_v)_{v\in V_n}$ is less than $d$, and Theorem 2.1 in
\citet{Janson04} shows that for any $t>0$ and $p\in[0,1]$,
\[
\mu_{n,p}\bigl(N<\EE_p(N)-t\bigr)\leq e^{-{2t^2}/({nd})} .
\]
On the other hand, by the bounded-degree assumption,
$\EE_\gamma(N) \ge(1-\gamma)^d n$, and therefore,
for any $c>1-(1-\gamma)^d$,
\[
\mu_{n,\gamma}\bigl(L_n^{(1)}\geq cn\bigr)\leq\mu_{n,\gamma}\bigl(N< (1-c)n\bigr)\leq
e^{-2n((1-\gamma)^d-(1-c))^2/d} .
\]
Choose $c$ such that
\[
1-(1-\gamma)^d+\sqrt{\frac{d\log({1}/{\eps})}{2n}} \le c < 1
\]
(which is possible for $n$ large enough), and get the desired
inequality.\vadjust{\goodbreak}
\end{pf}

The last piece we need for the proof is the fact that
the mean grows at least linearly in the super-critical phase.
This can be proved using the expansion property as follows:\vspace*{-3pt}
\begin{lemm}
\label{lemm:grandederivee}
Let $\alpha\in\ ]0,1/2[$ and $c\in\ ]0,1[$. There is a positive constant~$C'$,
depending only on $\alpha$, $c$, $b$ and $d$, such that for
$n$ large enough, and for every $p\in
[p_{n,\alpha}(c),p_{n,1-\alpha}(c)]$,
\[
\frac{d\EE_p(L_n^{(1)})}{dp}\geq C'n .\vspace*{-3pt}
\]
\end{lemm}

\begin{pf}
Let us fix $0<c<c_2<1$. Writing $f(x)=L_n^{(1)}(x)$ and
using Russo's lemma (\ref{eq:Russo}),
\begin{eqnarray*}
\frac{d\EE_p(L_n^{(1)})}{dp}&=&\frac{1}{p(1-p)}\sum_{e\in E_n}\EE
_p\bigl[\bigl(f(X)-f\bigl(X^{(e)}\bigr)\bigr)_-\bigr] \\[-2pt]
&\geq&\frac{1}{(1-p)}\EE_p\bigl(\bigl|\partial_E\bigl(\mathcal{C}^{(1)}\bigr)\bigr|\bigr) .
\end{eqnarray*}
By the expansion property,
\[
\bigl|\partial_E\bigl(\mathcal{C}^{(1)}\bigr)\bigr| \ge
b\bigl( L_n^{(1)}\II_{L_n^{(1)}\leq n/2}+\bigl(n-L_n^{(1)}\bigr)\II
_{L_n^{(1)}> n/2} \bigr),
\]
and therefore, for any $p\in[p_{n,\alpha}(c),p_{n,1-\alpha}(c)]$,
\begin{eqnarray*}
\frac{d\EE_p(L_n^{(1)})}{dp}
&\geq&\frac{b}{(1-p)}\EE\bigl[ L_n^{(1)}\II_{L_n^{(1)}\leq
n/2}+\bigl(n-L_n^{(1)}\bigr)\II_{L_n^{(1)}> n/2}\bigr]\\[-2pt]
&\geq& bn\min\{c,(1-c_2)\}\mu_p\bigl(L_n^{(1)}\in[cn,c_2n[\bigr) .
\end{eqnarray*}
Now, thanks to Proposition \ref{prop:criticalpaway01}, we know that
there is some $q_2(c)<1$ such that for $n$ large enough,
$p_{n,1-\alpha}(c)\leq q_2(c)$. Thus, applying
Lemma~\ref{lemm:proche1} with $\gamma=q_2(c)$ and $\eps=\alpha/2$, there is
some $c_2(c)\in\ ]c,1[$ such that, for $n$ large enough,
$p_{n,1-\alpha}(c)\leq q_2(c)\leq p_{n,\alpha/2}(c_2)$. Therefore
$\mu_p(L_n^{(1)}\geq cn)\geq\alpha$ and $\mu_p(L_n^{(1)}\geq
c_2n)\leq\alpha/2$ for any $p\in
[p_{n,\alpha}(c),p_{n,1-\alpha}(c)]$. This leads to
\[
\frac{d\EE_p(L_n^{(1)})}{dp}\geq bn\min\{c,(1-c_2)\} \alpha/2
\]
for any $p\in[p_{n,\alpha}(c),p_{n,1-\alpha}(c)]$.\vspace*{-3pt}
\end{pf}

Now we are ready to wrap up our argument.\vspace*{-3pt}
\begin{pf*}{Proof of Theorem \ref{theo:sharpthreshold}}
Let $c<1$ and $\alpha<1/2$ be fixed positive numbers.
We show that there exists a constant $K=K(\alpha,c,b,d)$, such
that if $\eps_n = K \log^{-1/3} n$, then
\[
p_{1/2}(c) - p_\alpha(c) \le\eps_n.
\]
The proof that $p_{1-\alpha}(c) - p_{1/2}(c) \le\eps_n$ is
completely similar.\vadjust{\goodbreak}

Lemma \ref{lemm:variancederivee} and the trivial bound
$L_n^{(1)} \le n$, implies that no matter how~$\epsilon_n$ is chosen,
\[
\int_{p_{1/2}(c)-\eps_n}^{p_{1/2}(c)-{3\eps_n}/{4}}\Var
_p\bigl(L_n^{(1)}\bigr) \,dp\leq
C \frac{n^2}{\log n},
\]
where $C=C(b,d)$.
Thus there is some $q_1\in
[p_{1/2}(c)-\eps_n,p_{1/2}(c)-3\eps_n/4]$ such
that
\[
\Var_{q_1}\bigl(L_n^{(1)}\bigr)\leq\frac{4Cn^2}{\eps_n\log n} .
\]
Similarly, one finds $q_2$ such that
\[
q_2\in
[p_{1/2}(c)-\eps_n/2,p_{1/2}(c)-\eps_n/4] \quad\mbox{and}\quad
\Var_{q_2}\bigl(L_n^{(1)}\bigr)\leq\frac{4C n^2}{\eps_n\log n} .
\]
Observe that it suffices to prove that $q_1 \le p_\alpha(c)$.
Note that
$q_1+\frac{\eps_n}{4}\leq q_2\leq p_{1/2}(c)$.
Thanks to Lemma \ref{lemm:grandederivee}, there is a
constant $C'$, depending only on $\alpha$, $c$, $b$ and $d$ such that
for $n$ large enough,
\[
\EE_{q_2}\bigl(L_n^{(1)}\bigr)- \EE_{q_1}\bigl(L_n^{(1)}\bigr)\geq C'n\frac{\eps_n}{4} .
\]
On the other hand, denote by
$M_p$ the median of $L_n^{(1)}$ under $\mu_p$ (we assume
that it is of the form $k+1/2$, with $k\in\NN$, which ensures its
uniqueness). $M_p$~is an increasing function of $p$ and therefore
\[
cn\geq M_{p_{1/2}(c)}-\tfrac{1}{2}\geq
M_{q_2}-\tfrac{1}{2} .
\]
By L\'{e}vy's inequality, the difference between the mean and median of any
random variable is bounded by its standard deviation and therefore
\[
\bigl|\EE_{q_2}\bigl(L_n^{(1)}\bigr)-M_{q_2}\bigr|
\leq\sqrt{\Var_{q_2}\bigl(L_n^{(1)}\bigr)}
\leq
n\sqrt{\frac{4C}{\eps_n\log n}} .
\]
Summarizing, we can write
\begin{eqnarray*}
&&\mu_{q_1}\bigl(L_n^{(1)}\geq cn\bigr) \\
&&\qquad= \mu_{q_1}\bigl(L_n^{(1)}-\EE_{q_1}\bigl(L_n^{(1)}\bigr)
\geq cn-\EE_{q_1}\bigl(L_n^{(1)}\bigr)\bigr) \\
&&\qquad\leq\mu_{q_1}\biggl(L_n^{(1)}-\EE_{q_1}\bigl(L_n^{(1)}\bigr)\geq
M_{q_2}-\frac{1}{2}-\EE_{q_1}\bigl(L_n^{(1)}\bigr)\biggr)\\
&&\qquad= \mu_{q_1}\biggl(L_n^{(1)}-\EE_{q_1}\bigl(L_n^{(1)}\bigr)\geq
M_{q_2}-\EE_{q_2}\bigl(L_n^{(1)}\bigr)+\EE_{q_2}\bigl(L_n^{(1)}\bigr)-\EE
_{q_1}\bigl(L_n^{(1)}\bigr)-\frac{1}{2}\biggr)
\\
&&\qquad\leq
\mu_{q_1}\Biggl(L_n^{(1)}-\EE_{q_1}\bigl(L_n^{(1)}\bigr)\geq C'n\frac{\eps
_n}{4}-\frac{1}{2}-n\sqrt{\frac{4C}{\eps_n\log
n}}\Biggr) .
\end{eqnarray*}
Now we choose $K=(256C/({C'}^2\alpha))^{1/3}$, so
\[
\eps_n=\biggl(\frac{256C}{{C'}^2\alpha\log n}\biggr)^{1/3} .
\]
Clearly, for $n$ sufficiently large,
\[
C'n\frac{\eps_n}{4}-\frac{1}{2}-n\sqrt{\frac{4C}{\eps_n\log
n}}\geq C'n\frac{\eps_n}{8} .
\]
Thus, Chebyshev's inequality implies
\begin{eqnarray*}
\mu_{q_1}\bigl(L_n^{(1)}\geq cn\bigr)
&\leq&
\frac{\Var_{q_1}(L_n^{(1)})}{(C'n\eps_n/8)^2} \\
&\leq&
\frac{4Cn^2}{\eps_n\log n} \frac{64}{{C'}^2n^2\eps_n^2} \\
&= &\alpha.
\end{eqnarray*}
This implies that $q_1\leq p_\alpha(c)$, as desired.
\end{pf*}

The bound on the size of the threshold width in Theorem
\ref{theo:sharpthreshold} is quite likely not to be tight. Indeed,
one would rather be inclined to compare it to what happens in the
Erd\H{o}s--R\'{e}nyi random graph $\mathcal{G}(n,(1+\eps)/n)$. In this
case, the mean grows also linearly, and the fluctuations of the giant
component are approximately Gaussian,
with a variance of order $\Theta(n)$, which implies that the
threshold width is of order $\Theta(1/\sqrt{n})$.
A similar tight threshold seems to
hold in random $d$-regular graphs as well [see Theorem 3 in
\citet{Pittelregular08}].
Thus, it is natural to conjecture that the threshold is much smaller
in our setting as well. Note that if the underlying graph $G_n$ is
transitive,
then \citet{FriedgutKalai} implies that the threshold width is at most
$O(1/\log n)$, so our quantitative bound $O(1/\log^{1/3} n)$
seems to be very weak.\looseness=1

Finally, let us emphasize an open problem about ``exponential
decay'' of the probabilities of abnormally small or abnormally large
size of the giant cluster. Proposition \ref{prop:criticalpaway01}
implies that for any $c$ in $]0,1[$, for
$p$ close enough to~1, the probability of not having a component of
size at least $cn$ decays as quickly as
$e^{-\alpha n}$ for some $\alpha>0$. An open
problem is to find out whether this is the case as soon as
$p>p_n^*(c)+\eps$ for
some fixed $\eps>0$. A similar question can be posed on the left of
the threshold: at which speed does the probability of having a
component of
size at least $cn$ decay to zero when $p<p_n^*(c)-\eps$?\looseness=1

\section*{Acknowledgment}
We thank Gady Kozma for useful discussions.

\printaddresses

\end{document}